\documentclass[12pt,reqno]{amsart}

\usepackage[arrow,matrix,curve]{xy}

\usepackage[dvips]{graphicx} 

\usepackage{amssymb, latexsym, amsmath, amscd, array, hyperref,
  epigraph
}

\usepackage{breakurl}   

\theoremstyle{definition}

\numberwithin{equation}{section}

\title{Of pashas, popes, and indivisibles}

\author{Mikhail G. Katz} \address{Department of Mathematics, Bar Ilan
  University, Ramat Gan 5290002 Israel} \email{katzmik@math.biu.ac.il}

\author{David Sherry} \address{Department of Philosophy, Northern
  Arizona University, Flagstaff, AZ
  86011,US}\email{David.Sherry@nau.edu}
 
\author{Monica Ugaglia}
\address{Liceo Scientifico Galileo Galilei, Macerata, Italy}
\email{monica.ugaglia@gmail.com}

\begin{document}

\thispagestyle{empty}




\begin{abstract}
The studies of Bonaventura Cavalieri's indivisibles by Giusti,
Andersen, Mancosu and others provide a comprehensive picture of
Cavalieri's mathematics, as well as of the mathematical objections to
it as formulated by Paul Guldin and other critics.  An issue that has
been studied in less detail concerns the theological underpinnings of
the contemporary debate over indivisibles, its historical roots, the
geopolitical situation at the time, and its relation to the ultimate
suppression of Cavalieri's religious order.  We analyze sources from
the 17th through 21st centuries to investigate such a relation.
\end{abstract}

\maketitle

\tableofcontents


\epigraph{As regards the opinion on \emph{quantity made up of indi-
  visibles}, I have already written to the Provinces many times that
  it is in no way approved by me and up to now I have allowed nobody
  to propose it or defend it. \\ Mutio Vitelleschi, S.J.}

\epigraph{We consider that this proposition not only stands in
  opposition to the doctrine of Aristotle that is commonly assumed,
  but also cannot be proven in itself. \\ Jakob Bidermann, S.J.}

\epigraph{However, the stain of this doctrine was subsequently
  revealed, disturbing us concerning the mysteries of the Eucharist,
  as we are taught by the Church.
\\ Sforza Pallavicino, S.J.}

\section{Introduction}

Studies of Bonaventura Cavalieri's indivisibles by Giusti \cite{Gi80},
Andersen \cite{An85}, Mancosu \cite{Ma96}, Delgado \cite{De17} and
others provide a comprehensive picture of Cavalieri's mathematics, as
well as of the mathematical objections to it as formulated by Paul
Guldin and other critics.  An issue that has been studied less
systematically concerns the theological underpinnings of the
contemporary debate over indivisibles, its historical roots, the
geopolitical situation at the time, and its relation to the ultimate
suppression of Cavalieri's religious order.  We analyze sources from
the 17th through 21st centuries to investigate such a relation.

An attentive reader of a recent article on Cavalieri will perhaps be
puzzled by the following passage:
\begin{enumerate}\item[]
In 1668, Pope Clement IX suppressed the Jesuats, probably at the
request of the Republic of Venice to finance the fight against the
Ottoman Turks for the rule of the island of Crete.  [Amir] Alexander
{\ldots}\;has recently {\ldots}~argued that hostility of the powerful
Jesuits to Cavalieri's mathematics and everything it represented
played a role in this suppression.%
\footnote{Gonz\'alez et al.~\cite[p.\;26]{Go18}.}
\end{enumerate}
The jesuates were founded in 1360; Bonaventura Cavalieri (1598--1647)
was a jesuate, as was his student Stefano Degli Angeli (1623--1697).
The jesuits were founded in 1540.  In the decades between 1630 and
1668, opposing views on indivisibles were recognizably (though not
entirely) aligned with these two orders.

Mathematically inclined readers of this journal may have learned about
Cavalieri's principle%
\footnote{\label{f2}If two plane domains have the same height and the
  same cross-sectional length at every point along that height, then
  they have the same area.  A similar relation holds between volumes
  of solids and their cross-sectional areas.  Using a method
  resembling Cavalieri's principle, Archimedes was able to find the
  volume of a sphere given the volumes of a cone and cylinder in his
  work \emph{The Method of Mechanical Theorems}.}
in their undergraduate calculus classes.  What did Cavalieri's
mathematics represent that would provoke the hostility of the jesuits?
Why would the republic of Venice make an extraordinary request to
suppress a religious order?  Doesn't the practice of suppressing such
orders so as to finance military campaigns seem unusual?  Wherefore
would a pope accede to such a request?

We will report on what we have been able to discover concerning this
case of science operating in, and influenced by, a historical and
theological context going back to the 13th century.  We will first
seek to document jesuit opposition to Cavalieri's indivisibles in
Sections~\ref{s2} and \ref{s9}, and explore its animus in
Section~\ref{s4}.  Section~\ref{s10} deals with the contemporary
geopolitical situation, the 1668 suppression of the jesuates, and the
1669 Ottoman takeover of Candia.  Sections~\ref{s6} and \ref{s7} deal
with the aftermath.  Section~\ref{s8} contains suggestions for further
research.

\section{Cavalieri's books}
\label{s2}

In this section, we review Cavalieri's books on the geometry of
indivisibles and their influence on his contemporaries.

Bonaventura Cavalieri was encouraged to study mathematics by his
teacher Benedetto Castelli (a benedictine).  They met at the jesuate
convent San Girolamo
%
%
in Pisa in 1616.%
\footnote{Giusti \cite[p.\;2]{Gi80}, Drake \cite[p.\,115]{Dr81},
Dufner \cite[pp.\;351--352]{Du75}, Feingold \cite[p.\,12]{Fe03}.}
%
%

Galileo composed a treatise on continuous quantity (now lost) as early
as 1609.%
\footnote{Drake \cite[p.\;245]{Dr81b}.}
Furthermore, 
\begin{enumerate}\item[]
Cavalieri, who took his start from Galileo's analysis, importuned him
to publish that work in order that Cavalieri might proceed with the
publication of his own \emph{Geometry of Indivisibles}.%
\footnote{Ibid.}
\end{enumerate}

\subsection{New approach}
\label{s21b}

The earliest indication we have of Cavalieri's interest in a new
approach to calculating areas and volumes occurs in one of his first
letters to Galileo, dated 15 december 1621 \cite{Ca21}.
A substantial part of the contents of Cavalieri's book can already be
found in letters from 1622--23.%
\footnote{Jullien \cite[p.\;89]{Ju15b}.}
In 1623, in a letter now lost, Galileo expressed reservations
concerning Cavalieri's approach.%
\footnote{Delgado \cite[p.\;239]{De17}.}
(Galileo's reservations point in a different direction from that
implied by Alexander's claim that ``Letters from Cavalieri in the
following months suggest that Galileo, at the very least, encouraged
Cavalieri to continue his investigations.'')%
\footnote{Alexander \cite[p.\;87]{Al15}.  In fact, when Castelli's
  position at Pisa became available in 1626, Galileo recommended not
  Cavalieri but rather Niccol\`o Aggiunti to fill the position; see
  Giusti \cite[p.\,13]{Gi80}.
}
%
%
Manuscript copies of the book were circulating no later than 1627.%
\footnote{Festa \cite[p.\;310, note~(7)]{Fe92a}, Arrighi \cite{Ar73},
  Andersen \cite[p.\;296]{An85}.}
Cavalieri used the term \emph{indivisible} in reference to his
approach for the first time in 1627.%
\footnote{Delgado \cite[p.\;31]{De17}.}
The book \emph{Geometria indivisibilibus continuorum nova quadam
  ratione promota}~\cite{Ca35} appeared in 1635.  De Gandt observed
that
\begin{enumerate}\item[]
[Cavalieri] tried, as very few did, to cast his new instruments in the
mold of classical Euclidean exposition.  It is striking to see how
``Euclidean'' Cavalieri is, compared to other creators in mathematics
at the same time.%
\footnote{De Gandt  \cite[p.\,160]{De92}.}
\end{enumerate}
As noted by de Gandt, Cavalieri devoted a lot of attention to a
mathematical justification of his approach.%
\footnote{Developed by a devout jesuate, Cavalieri's indivisibles were
  neither ``materialistic in tone'' \cite[p.\;394]{Al18} nor ``riddled
  with paradoxes'' (ibid.), a mathematical claim for which historian
  Alexander provides no convincing evidence; see further in
  Sections~\ref{s44} and \ref{s46}.}

The book provoked several responses by jesuit scholars.  High-placed
jesuit scholars like Paul Guldin published tracts highlighting alleged
inconsistencies of the indivisibles; see Section~\ref{s9}.

\subsection{\emph{Exercitationes} (1647) and philosophy}
\label{s21}

Cavalieri sought to demonstrate that alleged paradoxes attributed to
his indivisibles stem from their unrestricted use (without respecting
the conditions he specified).%
\footnote{E.g., in his 5 april 1644 letter to Torricelli; see Delgado
  \cite[pp.\,149--151]{De17}.}
Cavalieri published responses to some of the criticisms in his 1647
book \emph{Exercitationes geometricae sex}.  Here he wrote:


\begin{enumerate}\item[]
You know the arguments I used to prove the first Method of the
\emph{Indivisibili},
and I think to defend it from criticisms.  But if this is not enough,
you have my other method, free of the notion of infinite indivisibles
taken together.
And if you are not satisfied,
consider the arguments I used to validate the use of the indivisibles
following the style of Archimedes, mainly in Exerc.\;II.\, In this way
I think that you will be free from any doubt.%
\footnote{Cavalieri, Exercitatio Tertia, \cite[p.\;241]{Ca47};
  translation ours.  A Russian translation can be found in
  \cite[pp.\;61--62]{Ca40}.}
\end{enumerate}
Alluding to the distinction between \emph{discovering} and
\emph{proving} widely discussed at the time, Cavalieri continues:
\begin{enumerate}\item[]
Even [Paul] Guldin admitted that they are very useful for discovering
things (p.\;331 lib.\;IV).  You have to decide if they can be accepted
also for proving, and in what form: as I used them or with reference
to the Archimedean method, by which my discoveries have been
validated.  Finally, use it as you prefer: in this kind of discussion,
more philosophical than geometrical, I do not want to lose any time.
The Geometers will judge{\ldots}%
\footnote{Ibid.}
\end{enumerate}

The thrust of Cavalieri's remark about philosophy is that the
discussion concerning which approach is preferable -- indivisibles or
Archime\-dean exhaustion -- is a philosophical rather than a
geometrical issue.  The discussion not being a mathematical one,
Cavalieri has no interest in entering it.  Notably, he does not assert
that rigor is the concern of philosophy rather than of geometry.  The
term \emph{rigor} does not occur in the passage.

Cavalieri's point concerns the discovery of a new useful approach
exploiting indivisibles.  He employs it to obtain new results.  If he
is requested to \emph{justify} his findings using traditional methods,
namely Archimedean exhaustion, he is able to do so.  However, he
discovered them in a different way.  Is his way acceptable?  Cavalieri
the geometer prefers to leave such questions to the philosophers.

\subsection{The Boyer--Kline `affair'}
\label{s23}

Inexplicably, Cavalieri's comment on philosophy has been distorted out
of recognition by numerous commentators starting with Carl Boyer.
Boyer claimed the following:

\begin{enumerate}\item[]
Rigor, he [Cavalieri] said, was the \emph{affair} of philosophy rather
than geometry.$^{106}$ \cite[p.\,123]{Bo59} (emphasis added)
\end{enumerate}
Boyer's footnote 106: ``\emph{Exercitationes geometricae sex},
p.\;241,'' attached to this claim, leaves no doubt that he is
referring to the passage we quoted from page 241 of Cavalieri's
\emph{Exercitationes} in Section~\ref{s21}.%
\footnote{Boyer's dismissive comment is typical of some historians'
  tendency to attribute to Cavalieri the role of an early purveyor of
  unrigorous mathematics, and to attribute more validity than is
  warranted to the criticisms by Cavalieri's jesuit opponents.}
Boyer's claim was picked up by Morris Kline, who paraphrased Boyer by
adding two pairs of quotation marks, and claimed the following
concerning Cavalieri:

\begin{enumerate}\item[]
``Rigor,'' he said, ``is the concern of philosophy not of geometry.''%
\footnote{Kline \cite[p.\,133]{Kl80}.}
\end{enumerate}
Kline provided no source for his claim (which would be literally true
if the third-person pronoun were taken to refer to Boyer).  The
Boyer--Kline distortion of Cavalieri has been since copied by numerous
sources in multiple languages.%
\footnote{Thus, Brisson and Ofman have two versions of their erudite
article: the archive version \cite{Br20} with 140 footnotes, and the
shorter published version \cite{Br22}.  The archive version not only
reproduces the Boyer--Kline distortion, but actually claims that Boyer
is ``quoting Cavalieri'' \cite[note\;88, p.\,15]{Br20} (the claim is
not found in the published version).  O'Connor and Robertson furnish
the following purported Quotation: ``Rigour is the affair of
philosophy, not of mathematics'' \cite{Oc14}, but like Kline provide
no source.}
%
While lacking a source, the portrayal of Cavalieri as unrigorous
dovetails with Boyer's narrative of the conceptual development of
mathematics toward the predetermined pinnacle reached in the 1870s by
the ``great triumvirate.''%
\footnote{Boyer \cite[p.\;298]{Bo59}.}

We should note that Cavalieri at no time seeks to sum `infinitesimal
magnitudes' even though his language sometimes seems to point in that
direction: `All the lines\ldots', `All the surfaces\ldots' In
particular, the technique of discarding negligible (infinitesimal)
terms -- that would prove to be controversial in the history of the
calculus -- was never used by Cavalieri.%
\footnote{It is therefore inexact to claim, as Brisson and Ofman do in
  the archive version of their text, that in Cavalieri, ``Geometrical
  surfaces are infinite sums of parallel lines that can be treated, in
  some way, as finite sums'' \cite[p.\,15]{Br20} (this claim is
  similarly not found in the published version \cite{Br22}).}
Starting in 1641, Evangelista Torricelli will reconceptualize
indivisibles by assigning to them width.%
\footnote{See e.g., Delgado \cite[pp.\,15--16, 117--121]{De17}.}
When Torricelli evaluates areas of plane figures, he uses only the
lines corresponding to the partition points of an infinitesimal
partition, rather than `all the lines'.

\subsection{Adequality in Cavalieri and Fermat}
\label{s24}

Interestingly, Cavalieri used the term \emph{adequal} to describe the
relation between the collection of lines on the one hand and the space
they occupy on the other.%
\footnote{``magnitudinem, quae adequatur spatio ab eisdem lineis
  occupato, cum illi congruat'' Cavalieri \cite[p.\,17]{Ca35};
  cf.~Muntersbjorn \cite[p.\;235]{Mu00}.}
He did not mean it in the sense of \emph{approximately equal}.
Breger's claim that in Fermat, the term \emph{adequal} does not mean
``to be almost equal'' is endorsed by Muntersbjorn,%
\footnote{Muntersbjorn \cite[note\;5, p.\;250]{Mu00}.}
while she admits that Fermat used the term when performing quadratures
via a bunch of parallelograms.%
\footnote{Op.\;cit., p.\;243.}
Such a procedure (\emph{not} used by Cavalieri) necessarily involves
discarding negligible terms.  Itard noted the use of the term by
Cavalieri, observed that Fermat used it in various senses,%
\footnote{Itard \cite[p.\,117]{It75}.}
and clarified such senses:
\begin{enumerate}\item[]
Lorsque l'on compare deux grandeurs de m\^eme esp\`ece on peut le
faire soit en \'etudiant leur diff\'erence, soit leur rapport.  On
pourra les consid\'erer comme voisines si l'on trouve leur
diff\'erence voisine de z\'ero, ou leur rapport voisin de $1$ {\ldots}
Fermat donne ainsi au concept le sens d'infiniment petits
\'equivalents, pour utiliser un langage leibnizien.%
\footnote{Itard \cite[pp.\;339--340]{It74}.}
\end{enumerate}
Accordingly, Fermat, unlike Cavalieri, did use the term in the sense
of a relation of infinite proximity.%
\footnote{Briefly, the original term is actually Diophantus'
  \emph{parisotes}, which does refer to an approximate equality.
  Bachet, when translating Diophantus, performed a semantic calque,
  rendering it by a Latin term that corresponds to \emph{adequality}.
  Following Bachet, Fermat used the term in the same sense to refer to
  his relation, which is arguably a relation of infinite proximity.
  See further in \cite{13e}, \cite{18d},
  \cite[Sections~3.2--3.5]{20e}, and note~\ref{f52}.}

In sum, Cavalieri's books introduced a novel approach to problems that
have traditionally been approached via the method of exhaustion.
Torricelli carried the method further.  Subsequent geometers were
familiar with Cavalieri's method mainly through Torricelli's work.

\section{Four jesuits against the indivisibles}
\label{s9}

The first half of the 17th century witnessed various attempts to teach
indivisibles in jesuit colleges.  Among the jesuits of the period,
conflicting attitudes existed toward the new techniques.  The
rank-and-file were attempting to use them in their teaching, and sent
repeated requests for authorisation, from their various locations in
the periphery, to the central authorities.  Palmerino writes:
\begin{enumerate}\item[]
[I]f one were to study systematically the many \emph{Censurae
  opinionum} that preceded the \emph{Ordinatio} of 1651 throughout the
first half of the century, one could reconstruct in quite some detail
the chronological and geographical spread of the new sciences all over
Europe and their surreptitious entrance through the back door of the
Jesuit educational institutions.%
\footnote{Palmerino \cite[p.\,188]{Pa03}.}
\end{enumerate}
See further on the \emph{Ordinatio} in Section~\ref{s47}.  Gr\'egoire
de Saint Vincent was among the jesuit mathematicians who worked with
indivisibles.  Among the jesuits who sought to discredit indivisibles
mathematically were Paul Guldin, Andr\'e Tacquet, Antoine de
Lalouv\`ere, and Mario Bettini, dealt with respectively in
Sections~\ref{s31}, \ref{s32}, \ref{s33}, and \ref{s34}.

\subsection{Guldin}
\label{s31}

Paul Guldin's critique of Cavalieri's indivisibles appeared in the
fourth (and last) book of his \emph{De centro gravitatis} (also called
\emph{Centrobaryca}), published in 1641 \cite{Gu41}.  Mancosu writes:
\begin{enumerate}\item[]
Guldin is taking Cavalieri to be composing the continuum out of
indivisibles, a position rejected by the Aristotelian orthodoxy as
having atomistic implications.%
\footnote{\label{f28}For an elaboration of Mancosu's comment on Guldin
  and atomism, see Section~\ref{s44}.}
\ldots{} Against Cavalieri's proposition that ``all the lines" and
``all the planes" are magnitudes -- they admit of ratios -- Guldin
argues that ``all the lines {\ldots}~of both figures are infinite; but
an infinite has no proportion or ratio to another infinite.''%
\footnote{Mancosu \cite[p.\;54]{Ma96}.  For Guldin's original Latin
  see \cite[note\;15, p.\;200]{Fe92b} or \cite[note~456,
    p.\;204]{De17}.}
\end{enumerate}
If, in computing areas, an infinite can have no proportion or ratio to
another infinite as Guldin claimed, then Cavalieri's principle%
\footnote{See note~\ref{f2}.}
certainly risks being unsound.  According to Cavalieri, a key
objection of Guldin's was that indivisibles only appear in the
continuum \emph{in potentia} (potentially) and not \emph{in actu},%
\footnote{Delgado \cite[p.\,163]{De17}.}
harking back to Aristotelian strictures on philosophical discourse
concerning infinity.

It does not seem to have been noticed by commentators that some of
Guldin's arguments are refutable even by 17th century standards.
Relying upon Euclid's Book~5, Definition~5 (known as Definition~4 in
modern editions), Guldin argues as follows:
\begin{enumerate}\item[]
{\ldots}~For the parts both of this figure and of another can be
multiplied, so that the small parts surpass the bigger ones, and
vice-versa, and the same applies to the lines that designate them, and
to what I said on the subject of lengths; and therefore, according to
the fifth definition of Book 5 of Euclid's Elements, they are said to
have a ratio to each other.%
\footnote{Guldin as translated by Radelet-de Grave
  \cite[p.\;80]{Ra15}.}
\end{enumerate}
Having invoked what he takes to be Euclid's authority, Guldin proceeds
to attack Cavalieri's method of indivisibles in the following terms:
\begin{enumerate}\item[]
And so, if we sum any accumulation either of lines, or of parts of
planes of a surface, they are finite and have a ratio to each other,
the lines obviously to the lines, and the parts to the parts, but not
the lines to the parts.  This is why, since there is not, and cannot
be, any accumulation of all the lines or of all the parts, neither can
there be any ratio between them, etc.%
\footnote{Ibid.}
\end{enumerate}
Euclid's \emph{definition} is closely related to the Archimedean
property.%
\footnote{See e.g., Giusti \cite[p.\;33 and note 17 there]{Gi80}.}
In his argument against Cavalieri's method, Guldin did not
relate to the fact that it is a \emph{definition} rather than an
\emph{axiom}.  Guldin's critique based on Book V Definition 5 is,
strictly speaking, not applicable to Cavalieri's method, because
Cavalieri was working with \mbox{codimension-1} indivisibles%
\footnote{See further in note~\ref{f26}.}  
rather than with infinitesimals \`a la Kepler and Leibniz (which
violate the Archimedean property when compared to $1$).  Guldin's
argument, which would in principle apply had Cavalieri dealt with
infinitesimals, is refutable.  Half a century later, Leibniz will make
it clear that his incomparable infinitesimals violate the
comparability notion put forward in this definition.%
\footnote{Leibniz \cite[p.\;322]{Le95b}.  For further references see
  note~\ref{f143}.}

Significantly, Guldin also accused Cavalieri of plagiarizing Kepler.%
\footnote{Mancosu \cite[p.\;51]{Ma96}.  In short, Guldin accused
  Cavalieri of copying nonsense from Kepler.  Radelet-de Grave
  mentions ``Guldin's skill in successfully criticizing Cavalieri's
  method whilst simultaneously praising the one from which his
  inspiration was drawn, namely Kepler's'' \cite[p.\;77]{Ra15}.  She
  also notes that ``Giusti was surprised, quite rightly, that
  Cavalieri did not take the opportunity to counter that Guldin's rule
  was already to be found in Pappus'' \cite[p.\;76]{Ra15}.  Cavalieri
  did point out the similarity of Guldin's rule to Kepler's; see
  \cite[p.\;217]{De17}.}
Given the dubiousness of some of Guldin's mathematical criticisms, one
can wonder what his true motivations were and whether there may have
been a hidden doctrinal agenda.  The fact that he would pursue
non-mathematical criticisms (such as allegations of plagiarism) lends
credence to the hypothesis that his opposition to Cavalieri's
indivisibles was not purely mathematical.  As Festa, Mancosu, and
others have pointed out, a well-known allusion by Guldin to the
existence of extra-mathematical issues may have involved a doctrinal
issue; see Section~\ref{itemfour}.  See further in Section~\ref{s4} on
the connection to atomism.

Cavalieri observed that, in spite of Guldin striking the posture of a
defender of mathematical certainty, the justification Guldin provided
for his own rule (for bodies of revolution)%
\footnote{See Delgado \cite[p.\,174, note~368]{De17}.}
is in fact insufficient.  Torricelli agreed with Cavalieri's
assessment.%
\footnote{Delgado \cite[pp.\,174--175]{De17}.}
Cavalieri went on to point out that Antonio (Giannantonio) Rocca
(1607--1656) furnished a better-grounded justification for the rule in
terms of{\ldots} indivisibles,%
\footnote{Delgado \cite[pp.\,178--179]{De17}; Giusti \cite[p.\;50,
    note 40]{Gi80}.}
while Cavalieri himself provided a generalisation for solids.%
\footnote{Delgado \cite[p.\,180]{De17}.}
Of interest also is an error, discovered by Torricelli, on page 127 of
the first book of Guldin's work~\cite{Gu41}.%
\footnote{Giusti \cite[p.\;75]{Gi80}; Delgado
  \cite[pp.\,184--187]{De17}.}

When Cavalieri first became aware of Guldin's criticism in 1642, he
made the following visionary comment:
\begin{enumerate}\item[]
I do not mind if this Father has taken upon himself to refute this
method of mine on indivisibles, because, if I am in error I will see
the truth, but if he is the one who is mistaken, at least he will do
this favour to my Geometria, which is that some people, who otherwise
might never have seen it, will make some reflection upon it.%
\footnote{Cavalieri as translated by Giusti \cite[p.\;55]{Gi80}.}
\end{enumerate}

\subsection{Tacquet}
\label{s32}

Andr\'e Tacquet was generally less extreme than some of the other
jesuit opponents of indivisibles.  He tolerated their use in heuristic
arguments leading to the discovery of results that must then be
justified by traditional exhaustion arguments.  As noted by Bosmans,
\begin{enumerate}\item[]
[Tacquet] connaissait la \emph{G\'eom\'etrie des indivisibles} de
Cavalieri et appr\'eciait la f\'econdit\'e des m\'ethodes du
g\'eom\`etre italien.  Elles \'etaient cependant, d'apr\`es lui, un
simple moyen de recherche, sans force d\'emonstrative, \`a moins qu'on
ne les ramen\^at \`a la m\'ethode d'exhaustion d'Archim\`ede.%
\footnote{Bosmans \cite[p.\;71]{Bo27}.}
\end{enumerate}
Nonetheless, he warned in his \emph{Cylindricorum et annularium libri
  IV} (1651) that the idea of quantity composed of indivisibles makes
war upon geometry to such an extent that ``if it is not to destroy it,
it must itself be destroyed''%
\footnote{Festa \cite[p.\;205]{Fe92b}; see note 26 there for the
  original.  Cf.~Descotes \cite[p.\;255]{De15}.}
(i.e., the method of indivisibles must itself be destroyed if it is
not to destroy geometry).

Tacquet's terminology of ``quantity composed of indivisibles'' echoes
the language of the 1632 ban, suggesting that the ban applied to
mathematical indivisibles, as well; see further in Section~\ref{s44}.

Giusti provides the following apt summary:
\begin{enumerate}\item[]
Tacquet rejects as non-geometrical the theory of indivisibles, which
he calls the method of proof \emph{per heterogenea}.  However, he
proves many of his own theorems both by indivisibles and by classical
Archimedean methods.%
\footnote{Giusti \cite[note 49, p.\,21]{Gi80}.}
\end{enumerate}

\subsection{Lalouv\`ere}
\label{s33}

A declared enemy of indivisibles, jesuit Antoine de Lalouv\`ere
unconsciously used them in his own work, as noted by Descotes:
\begin{enumerate}\item[]
What is quite striking in this procedure is that Lalouv\`ere, who is a
declared enemy of indivisibles, uses without even realizing it an
equivalent, and even at heart an identical, process to Dettonville's,%
\footnote{Dettonville is a pen name for Pascal.}
given that all the surfaces are arranged on the same line~$AC$,
divided into small portions. This amounts to summing these magnitudes,
taken all along the line~$AC$, while disregarding the problem posed by
the divisions of~$AC$, and consequently, the presence of indivisibles.
{\ldots}
%
%
At the very moment when he attacks the indivisibles, Father
Lalouv\`ere unconsciously uses them.%
\footnote{Descotes \cite [p.\;273]{De15}.  Giusti similarly points out
  that Lalouv\`ere ``did not reject completely the method [of
    indivisibles] but tried to prove it by exhaustion, a proof much in
  the spirit of [Cavalieri's] \emph{Geometria}, and of the
  \emph{Exercitatio} II'' \cite[p.\;45]{Gi80}.}
\end{enumerate}
Lalouv\`ere is one of the jesuit mathematicians analyzed in a
comprehensive study of the period 1540--1640 by Antonella Romano.  She
writes:
\begin{enumerate}\item[]
Un premier constat regarde les j\'esuites confront\'es \`a la censure:
sur tous les cas du Fondo Gesuitico qui concernent la France, deux
professeurs seulement appartiennent \`a la liste \'etablie dans le
cadre de cet ouvrage, B. Labarthe et V. L\'eotaud.  Si d'autres
math\'ematiciens y apparaissent, c'est tout aussi exceptionnellement,
et au titre de censeur, comme Antoine Lalouv\`ere.$^{114}$
\cite[p.\;512]{Ro99}.
\end{enumerate}
Romano's footnote\,114 there reads: ``C'est lui qui porte un jugement
n\'egatif sur l'ouvrage de B. Labarthe.''  Thus, Lalouv\`ere in his
capacity of censor sank at least one mathematical book, namely the
\emph{Hermetis Mathematici praeludium}%
\footnote{Romano \cite[p.\;515]{Ro99}.}
of his fellow jesuit Bartholom\'e Labarthe, in~1662.%
\footnote{The fact that jesuits Nicolas du Port Guichart and Jean-Paul
  M\'edaille
%
%
considered Labarthe's book to be ``correct'' \cite[p.\;515]{Ro99}
suggests that Lalouv\`ere's opinion may have been controversial.}
On the relation of Lalouv\`ere and Fermat see the note.%
\footnote{\label{f52}The presence of an opinionated jesuit
  anti-indivisibilist at Toulouse may help explain Fermat's reticence
  when discussing the foundations of his method of adequality and the
  nature of the increment~$E$; see \cite{13e}, \cite{18d}.
  Lalouv\`ere was at Toulouse during the periods 1632--35, 1643--54,
  and 1659--64 \cite[pp.\;584--585]{Ro99} (and taught logic in
  1632--33 and physics in 1633--34 \cite[p.\;548]{Ro99}).  Fermat
  deposited his method of adequality with d'Espagnet already in the
  late 1620s.  The fact that Fermat did not communicate it to Mersenne
  in Paris until 1636 \cite[p.\;53]{Ma94} - after Lalouv\`ere's
  departure from Toulouse the previous year, 1635 - may have been a
  coincidence.}

\subsection{Bettini}
\label{s34}

To jesuit Mario Bettini, Cavalieri's indivisibles were not merely
`hallucinations' but glib philosophizing that makes geometrical
theorems useless:
\begin{enumerate}\item[]
{\ldots}\;being pressed, I respond to the counterfeit philosophizing
about geometrical figures by indivisibles.  Far, far be it from me to
wish to make my geometrical theorems useless, lacking demonstrations
of truth.  Which would be to compare {\ldots}\;figures and
philosophize about them by indivisibles.%
\footnote{Bettini as translated by Alexander \cite[p.\,159]{Al15}.}
%
%
\end{enumerate}
The criticism appeared in Bettini's 1648 \emph{Aerarium philosophiae
  mathematicae} \cite{Be48}.  A sharp-tongued rebuttal by Degli Angeli
in \cite{De58} is analyzed by Alexander.%
\footnote{Alexander \cite[p.\,168]{Al15}.}
Degli Angeli's response to Tacquet in \cite{De59} is analyzed by
Festa.%
\footnote{Festa  \cite[p.\;205]{Fe92b}.}

In sum, some leading jesuit mathematicians attacked indivisibles
ostensibly on mainly mathematical grounds (though there is an allusion
in Guldin to broader issues; see Section~\ref{itemfour}).  In
Section~\ref{s4}, we examine the related actions by their doctrinal
superiors in the order.

\section{Bans on quantity made up of indivisibles}
\label{s4}

Parallel to the efforts sketched in Section~\ref{s9} by jesuit
scholars such as Bettini, Guldin, Lalouv\`ere, and Tacquet to
discredit indivisibles mathematically, there was a series of bans
against the teaching of indivisibles, issued by their doctrinal
superiors in the order.  The bans are catalogued by Festa \cite{Fe92b},
Palmerino~\cite{Pa03}, and Alexander~\cite{Al15}.

\subsection{Physical \emph{vs} mathematical indivisibles}
\label{s41}

The jesuit bans applied equally to physical atomism and mathematical
indivisibles.  On 10~august 1632, the Revisors General led by Jakob
Bidermann formulated and then banned the following proposition:
\begin{enumerate}\item[]
%
%
A permanent continuum can be constituted by physical indivisibles
alone, or atomic [i.e., non-divisible] particles, having
\emph{mathematical parts} that can be designated, even though the
mentioned particles are distinct from each other.
%
%
Time, too, [consists] of instants, and intense qualities consist of
indivisible grades alone.

%
%
We consider that this proposition not only stands in opposition to the
doctrine of Aristotle that is commonly assumed, but also cannot be
proven in itself; {\ldots}%
\footnote{ ``Continuum permanens potest constare ex solis
  indivisibilibus physicis, seu corpusculis atomis, habentibus partes
  mathematicas, in ipsis designabiles, etiamsi realiter dicta
  corpuscula inter se distinguantur.  Tempus quoque ex instantibus, \&
  qualitates intensae, ex solis gradis indivisibilibus constant.

Hanc propositionem arbitramur, non modo repugnare communi Aristotelis
doctrina, sed etiam secundum se esse improbabilem, etc.''
\cite[p.\;207]{Fe92b}.  For a French translation see
\cite[p.\,198]{Fe92b}.}
(emphasis on `mathematical parts' added)
\end{enumerate}
The explicit reference to Aristotelian doctrines is significant.
Aristotle, unlike Plato, held that mathematical objects were all
instantiated in physical objects, or that mathematical objects are
abstracted from physical objects.%
\footnote{E.g., on Aristotle's view, a sphere is abstracted from the
  moon, by abstracting the sensible property (sphericity, say) from
  its matter.  See further in Lear \cite{Le82}.}
Accordingly, Aristotelian philosophy of mathematics would couple
mathematical and physical indivisibles.  What did the Revisors General
mean exactly by `mathematical parts' of physical indivisibles?  Festa
speculates that Cavalieri's indivisibles may be alluded to in this
ban%
\footnote{Manuscript copies of Cavalieri's book were already in
  circulation in 1627--1629, thus several years prior to the ban; see
  Section~\ref{s21b}.}
\cite[p.\,199]{Fe92b}.  As noted by Mancosu,
\begin{enumerate}\item[]
[Guldin's] preface defines mathematics according to the Aristotelian
classification, as that part of philosophy lying between physics and
metaphysics.  Mathematics is the science that considers quantity
abstracted from sensible matter.%
\footnote{Mancosu \cite[p.\;56]{Ma96}.}
\end{enumerate}
Guldin's mathematical critique of indivisibles was examined in
Section~\ref{s31}.  Guldin, too, would likely have interpreted the ban
as applying equally to mathematical and physical indivisibles.%
\footnote{Guldin's book may contain an allusion to a doctrinal issue;
  see further in Section~\ref{itemfour}.}

Referring to the 1632 ban, Feingold notes:
\begin{enumerate}\item[]
Six months later, General [Mutio] Vitelleschi formulated his strong
opposition to mathematical atomism in a letter he dispatched to Ignace
Cappon in Dole: ``As regards the opinion on \emph{quantity made up of
  indivisibles}, I have already written to the Provinces many times
that it is in no way approved by me and up to now I have allowed
nobody to propose it or defend it.''%
\footnote{Feingold \cite[pp.\;28--29]{Fe03}; emphasis added.}
\end{enumerate}
Significantly, Vitelleschi speaks of \emph{quantity}, which Aristotle
abstracts from matter, confirming the coupling of mathematical and
physical indivisibles in the thinking of the General.

Given such an explicit linkage, it is difficult to agree with
Alexander's view that specifically \emph{physical} indivisibles were
the target of the bans which, as he claims, lacked grounds to object
to mathematical ones:
\begin{enumerate}\item[]
In as much as the technical details of the miracle of the Eucharist
mattered, they provided no grounds for objecting to a mathematical
(rather than physical) doctrine.%
\footnote{Alexander  \cite[p.\;393]{Al18}.}
\end{enumerate}
See further in Section~\ref{s45b} for an analysis of Alexander's
position.  For the precise canon involved see Section~\ref{f7}.

\subsection{Atomism \emph{vs} canon}
\label{f7}

What was behind the jesuit opposition to indivisibles?  Sforza
Pallavicino pinpointed the problem as a clash with (the catholic
interpretation of) the eucharist.  Festa recounts jesuit Pallavicino's
1647 recanting of his previous interest in physical atomism:
\begin{enumerate}\item[]
Cette doctrine, \'ecrit Sforza Pallavicino, ``flatteuse pour notre
imagination, {\ldots}~attira une foule de disciples.  Mais son
caract\`ere destructeur a \'et\'e d\'etect\'e; elle trouble ce que
l'\'Eglise nous enseigne sur les Myst\`eres de l'Eucharistie et ne
s'accorde pas suffisamment avec ce que le Concile du Latran affirme
sur la nature de l'\^ame humaine.''%
%
\footnote{Pallavicino as translated by Festa \cite[item~10]{Fe99};
    cf.\;\cite[p.\;203]{Fe92b}.}
\end{enumerate}

The endorsement of transubstantiation at the 16th century Council of
Trent, Session 13, Canon 2 was widely interpreted by catholic
theologians as endorsement of the Peripatetic theory of matter and
form sometimes referred to as hylomorphism (or hylemorphism),%
\footnote{As a 0-th approximation for the benefit of a reader not
  familiar with the circle of ideas related to hylomorphism, note
  that, for example, a clay pot can be viewed as having two basic
  components: (1) undifferentiated \emph{matter}, or material (clay),
  and (2) its \emph{form} or shape.  Such a perspective creates
  tensions with atomistic approaches that go against the grain of
  thinking in terms of undifferentiated matter.}
which they saw as opposed to atoms, closely related to indivisibles.%
\footnote{On transubstantiation and the Council of Trent, see
  Armogathe \cite[pp.\;28--32]{Ar77}, McCue \cite[p.\;419]{Mc68},
  Festa \cite[p.\,101]{Fe91}, and Section~\ref{s8}.  On Scotus
  concerning the incompatibility of atomism and the eucharist, see
  Section~\ref{s47}.}
Canon~2 sought to clarify rulings of the 4th Lateran Council three
centuries earlier; see Section~\ref{s47}.  L\"uthy and Nicoli speak of
atoms becoming a problem
\begin{enumerate}\item[]
in the Counter-Reformational context of tightening doctrinal rules,
and more specifically {\ldots}~the decision of the Council of Trent to
insist on transubstantiation as the only correct way of interpreting
the Eucharist, wherein transubstantiation was defined in the
terminology of substance and accident, of matter and substantial
forms.%
\footnote{L\"uthy and Nicoli \cite[p.\;7]{Lu23ch1}.}
\end{enumerate}

\subsection{Grassi and Inchofer \emph{vs} Galileo}
\label{s43}

In the case of Cavalieri's mentor Galileo, there are multiple
documents by jesuits -- specifically, Orazio Grassi and Melchior
Inchofer -- attacking his atomism with the weapon of eucharistic
theology.  Thus, Ferrone and Firpo write:
\begin{enumerate}\item[]
[I]n his \emph{Ratio ponderum}, [jesuit Orazio] Grassi was to include
the relationship between atomism and the Eucharist among the many
accusations he leveled at [Galileo's book] \emph{Saggiatore}.%
\footnote{Ferrone and Firpo \cite[p.\;505]{Fe86}.}
\end{enumerate}
The Galileo--Grassi debate and its theological underpinnings are
analyzed by Festa.%
\footnote{Festa \cite[items~7--8]{Fe99}.}
A detailed denunciation, document labeled EE291, came to light in
1999.%
\footnote{Artigas et al.~\cite{AMS}, Finocchiaro
  \cite[p.\;291]{Fi21}.}
The author of EE291 is thought to be the jesuit Inchofer based on
handwriting analysis.%
\footnote{Festa \cite[p.\;25]{Fe07}.}
The criticisms voiced in EE291 echo the accusations against Galileo
leveled by both Grassi and an earlier anonymous document labeled~G3.%
\footnote{The centrality of the atomist issue to the background
deliberations in the trial of Galileo was the thesis of Redondi's
book~\cite{Re87}.  Our argument is independent of Redondi's thesis.
Namely, the relevant issue for us is Grassi's critique of Galileo
based on the relation between atomism and the eucharist, a historical
fact acknowledged by Redondi's critics Ferrone and Firpo.}
An alternative interpretation of the jesuit opposition to
indivisibles, proposed by Alexander, is examined in
Section~\ref{s45b}.

\subsection{Galileo and Cavalieri}

Significant differences existed between Galileo and Cavalieri with
regard to indivisibles.%
\footnote{Jullien \cite[pp.\;95--97]{Ju15b} based on Giusti
  \cite[pp.\;40--41]{Gi80}.}
Galileo was initially reluctant to use them in geometry.  Cavalieri,
on the other hand, thought of indivisibles as a kind of `useful
fiction' (of course a mathematical one), eminently applicable in
geometry.  His letters to Galileo suggest that Cavalieri was
disappointed with Galileo's attitude.  The long hiatus between the
writing of Cavalieri's book (1627) and its publication in extended
form (1635) was apparently due to his efforts to influence Galileo's
opinion concerning geometric applications of indivisibles.

Delgado argues that Cavalieri's use of indivisibles to provide a new
solution to the problem of free fall in his 1632 \emph{Lo Specchio
  Ustorio} \cite{Ca32} succeeded in stimulating Galileo's interest in
Cavalieri's geometric method.%
\footnote{Delgado \cite[pp.\;45--46]{De17}.}
Delgado traces Galileo's evolution from a physical to a mathematical
atomism.%
\footnote{Delgado \cite[p.\;55]{De17}.}
Palmerino analyzes Cavalieri's influence on Galileo.%
\footnote{Palmerino \cite[p.\;307]{Pa00}; cf.~Delgado
  \cite[p.\;79]{De17}.}
It is interesting to note that in his first reaction to the 1638
\emph{Discorsi} where Galileo proves the law of fall, Cavalieri
``criticizes Galileo for not having emphasized that the indivisibles
have to be taken as equidistant.''%
\footnote{Damerow et al.~\cite[note 175, p.\;251]{Da04}.}

\subsection{The Festa--Mancosu thesis}
\label{s44}

Historians are divided with regard to the reasons behind the jesuit
opposition to indivisibles.  Mancosu writes:
\begin{enumerate}\item[]
To these open attacks [on Cavalieri's work] one should however add the
general hostility of the Jesuits to Cavalieri's indivisibles.  Festa
(1990, 1992) provides archival evidence to show that the teaching of
indivisibilist techniques in geometry, as well as the use of the
atomistic theory in physics, was forbidden in the Jesuit schools by
means of decrees, the first dated 1632, issued by the `revision'
fathers of the Collegium Romanum.  Festa argues that the hostility to
atoms and indivisibles was motivated by theological concerns about the
dogma of transubstantiation.%
\footnote{Mancosu \cite[note\;18, pp.\;219--220]{Ma96}.}
\end{enumerate}
Accordingly, the Festa--Mancosu thesis is that the teaching of
indivisibilist techniques in geometry was included in the prohibition
due to doctrinal concerns over transubstantiation.

\subsection{Euclid's geometry as a world-ordering force}
\label{s45b}

Meanwhile, Alexander proposes an alternative interpretation of the
jesuit opposition to mathematical indivisibles.  He argues that to the
jesuits starting with Clavius, Euclid's mathematics served as a
world-ordering force that (as they were convinced) could be used in
educating people and attracting them to catholicism:

\begin{enumerate}
\item
``The War Against Disorder: The Jesuits against the infinitely
  small.''%
\footnote{Alexander \cite[p.\,15]{Al15}.  By mentioning the ``war
  against disorder'' already in the heading of an entire section,
  Alexander seeks to emphasize his idea that the jesuits were out to
  save the world from chaos and disorder, and impose upon it
  mathematical order based on Euclid.}
\item
``And like each proof alone, geometry as a whole is universally and
  eternally true, \emph{ordering the world} and governing its
  structure everywhere and always.''%
\footnote{Alexander \cite[p.\;67]{Al15}; emphasis on ``ordering the
  world'' added.}
\item
``It was clear to Clavius that Euclid's method had succeeded in doing
  precisely what the Jesuits were struggling so hard to accomplish:
  imposing a true, eternal, and unchallengeable \emph{order} upon a
  seemingly chaotic reality.''%
\footnote{Ibid.; emphasis added.}
\item
``The Jesuits {\ldots}\;believed that thanks to the rational rigor of
  its method and the unshakeable certainty of its results, mathematics
  would play a key role in imposing \emph{order} on the chaos brought
  about by the Reformation, and reestablishing the authority of the
  Church hierarchy.''%
\footnote{Alexander \cite[p.\;394]{Al18}; emphasis added.}
\end{enumerate}
Accordingly, Alexander claims that jesuit opposition to mathematical
indivisibles was prompted by the allegedly paradoxical nature of the
latter, inconsistent with the Euclidean ideal they sought to
safeguard.  On occasion, Alexander suggests that indivisibles are not
merely paradoxical but inherently contradictory, as when he writes:
\begin{enumerate}
\item[(i)] ``Zeno's mind-benders prove extremely difficult to resolve,
  based as they are on the inherent contradictions posed by
  indivisibles.''%
\footnote{Alexander \cite[pp.\;9--10]{Al15}.}
\item[(ii)] ``Most damaging of all, whereas Euclidean geometry was
  rigorous, pure, and unassailably true, the new methods were riddled
  with paradoxes and contradictions, and as likely to lead one to
  error as to truth.''%
\footnote{Alexander \cite[p.\,120]{Al15}.}
\end{enumerate}
Alexander may have magnified such an aspect of indivisibles so as to
buttress his interpretation, analyzed in Section~\ref{s46}.

\subsection{Analysis of Alexander's interpretation}
\label{s46}

Alexander's interpretation of jesuit opposition to indivisibles is
clearly at odds with the Festa--Mancosu thesis (see
Section~\ref{s44}).  Alexander's interpretation can be challenged.  We
note the following six points.

\subsubsection{Galbraith's study}
\label{s462}

The study by Galbraith~\cite{Ga21} of the history of the jesuits (in
relation to mathematics) suggests that the catholic hierarchy was
ambivalent about the involvement of members in mathematics in general
(including the geometry of Euclid), and that many among the higher-ups
thought that such involvement detracted members from their primary
responsibilities.  Accordingly, the attitude toward mathematics as a
world-ordering force may have been mainly the attitude of Clavius
himself%
\footnote{As noted by Hellyer, ``Clavius' textbooks standardized
  mathematical instruction throughout the Order to a very high degree
  well into the seventeenth century'' \cite[note 12, p.\;323]{He96}.
  See also \cite{Cl02}.}
(see further in Section~\ref{cl}).  Galbraith's analysis specifically
undermines Alexander's claim that
\begin{enumerate}\item[]
By the late sixteenth century, mathematics had become one of the most
prestigious fields of study at the Collegio Romano and other Jesuit
schools.%
\footnote{Alexander \cite[p.\,120]{Al15}.}
\end{enumerate}

\subsubsection{Clavius}
\label{cl}
Significantly, Clavius was a proponent of viewing hornangles as
quantities (an example of a hornangle is the crevice between a
circular arc and the tangent ray at one of its endpoints).  A
well-known controversy opposed Clavius and Peletier on the subject of
hornangles (with Clavius in favor and Peletier opposed); see e.g.,
\cite{Ma90}, \cite{Ax18}.%
\footnote{An additional controversy concerning the status of
  mathematics opposed Clavius and Perera (Pereira, Pereyra); see
  Mancosu \cite[p.\,14]{Ma96}, Rossini~\cite{Ro22}.}
Whereas Alexander argues that the reason for the jesuits' opposition
is their desire to safeguard the world-ordering power of Euclid
against an alleged threat stemming from the non-Archimedean behavior
of infinitesimals, hornangles represent a non-Archimedean phenomenon:
a hornangle is smaller than a half, a third, a quarter, etc. of any
ordinary rectilinear angle.  Therefore Clavius's advocacy of
hornangles undermines Alexander's claim that, following Clavius, the
jesuits were duty bound to oppose Cavalieri's mathematics.%
\footnote{Knobloch writes concerning the Clavius--Peletier
  controversy: ``Si nous laissons de c\^ot\'e sa querelle inutile avec
  Peletier au sujet de l'angle de contingence, il se montre le plus
  souvent comme un critique subtil rep\'erant les constructions
  fausses, les d\'emonstrations erron\'ees ou les propositions mal
  formul\'ees'' \cite[p.\;337]{Kn88}.  We feel that analyzing
  non-Archimedean phenomena is not ``useless'', for both Clavius and
  Leibniz; see \cite[Section~3]{21g}.}
%


\subsubsection{Ordinatio on incompatibility}

Alexander discusses the 1651 jesuit \emph{Ordinatio},%
\footnote{Alexander \cite[pp.\,147--149, 151, 321]{Al15}.}
but not the fact that it explicitly mentions the incompatibility of
atomism with the eucharist; see further in Section~\ref{s47}.  This
aspect of the \emph{Ordinatio} points in a different direction from
the one argued by Alexander.

\subsubsection{Guldin's allusion}
\label{itemfour}
As noted by Alexander,
\begin{enumerate}\item[]
Guldin comes close to admitting that there are greater issues at stake
than the strictly mathematical ones, writing cryptically that ``I do
not think that the method [of indivisibles] should be rejected for
reasons that must be suppressed by never inopportune silence,'' but he
gives no explanation of what those ``reasons that must be suppressed''
could be.%
\footnote{Alexander  \cite[p.\,154]{Al15}.}
\end{enumerate}
Guldin's mathematical critique of indivisibles was examined in
Section~\ref{s31}.  What issue with indivisibles could Guldin be
alluding to in the passage above?  If, as Alexander claims, Guldin is
alluding to the world-ordering power of the timeless rigor of Euclid's
mathematics, why should such a comment be suppressed by a `never
inopportune silence'?  On the contrary, Guldin should have emphasized
it!  We don't know exactly what Guldin had in mind, but he was likely
referring to extra-mathematical issues.  It is more plausible to
interpret Guldin's quip as alluding to doctrinal problems with
indivisibles than to their posing a perceived threat to the
world-ordering power of Euclidean rigor.

\subsubsection{Bans as evidence}

In his book \cite{Al15}, Alexander presents 17th century jesuit bans
against indivisibles as evidence for his thesis that jesuits were
opposed to mathematical indivisibles because indivisibles were
paradoxical and contradictory, and therefore undermined Euclidean
geometry, which, according to Alexander, the jesuits viewed as a
world-organizing force.  If one assumes that the bans applied only to
physical indivisibles,%
\footnote{As claimed by Alexander; see end of Section~\ref{s41}.}
there is a basic incoherence in Alexander's approach: if the jesuit
mathematicians were mainly bothered by the contradictory nature of
\emph{mathematical} indivisibles, then their doctrinal superiors' bans
on \emph{physical} indivisibles constitute an unrelated issue that
provides no evidence for Alexander's thesis.

\subsubsection{An old problem}

Catholicism's problems with indivisibles were not a
Counter-Reformation phenomenon (as implied by Alexander), and in fact
date back to the 13th century; see Section~\ref{s47}.

\subsection{Atoms, shibboleths, Descartes}
\label{s45}

L\"uthy notes that
\begin{enumerate}\item[]
At the latest by 1650, the atom had become a fetish to one party and a
scandal to the other, and hence a shibboleth that defined one's
adherence to the so-called new philosophy.%
\footnote{L\"uthy \cite[p.\,114]{Lu23}.}
\end{enumerate}
Rossi presents a detailed study of ``the radical incompatibility
between atomist theses and the conclusions reached by the Council of
Trent on the sacrament of the Eucharist'' \cite{Ro98}, and analyzes
attempts by maverick jesuit scholars (Suarez, Pereira, Arriaga, and
Oviedo) ``to develop a natural philosophy alternative to the
Aristotelian one, capable of retaining its distance from impious
atomism, but, nonetheless, adopting certain central aspects concerning
belief in discontinuity and indivisibles'' (ibid.).

Ren\'e Descartes questioned the explanatory power of prime matter and
substantial forms (the two components of hylomorphism) in the
following terms:
\begin{enumerate}\item[]
[W]e do much better to understand what takes place in small bodies,
whose minuteness prevents us from perceiving them, by what we see
occurring in those that we do perceive $<$and thus explain everything
in nature, as I have tried to do in this treatise$>$. This is
preferable to explaining certain things by inventing all sorts of
novelties with no relation to those that are perceived $<$such as
prime matter, substantial forms, and all the whole range of qualities
which many are in the habit of assuming, any one of which is more
difficult to understand than all the things they are supposed to
explain$>$%
\footnote{Descartes as translated by Ariew in \cite[Article~201,
    p.\;269]{De00}.}
\end{enumerate}
While Descartes was an anti-atomist,%
\footnote{\label{f19}Thus, in article 202 of his \emph{Principles of
Philosophy}, Descartes distanced himself from the atomism of
Democritus in the following terms: ``202.  Que ces principes ne
s'accordent point mieux auec ceux de Democrite qu'auec ceux d'Aristote
ou des autres'' \cite[pp.\;62, 516]{De24}.  See also
\cite[p.\;3]{Lu23ch1}.}
corpuscularism was a common feature of the approaches of Kepler,
Galileo, and Descartes.%
\footnote{Berchman \cite[pp.\,199--200]{Be09}.}

In 1663, some of Descartes's writings were placed on the Index of
Prohibited Books.%
\footnote{Armogathe and Carraud \cite[p.\;67]{Ar04}.}
Stephanus Spinula was one of the two censors:
\begin{enumerate}\item[]
By the time the Congregation of the Holy Office gives [Spinula] the
task of judging Descartes's \emph{Principles} and the \emph{Passions
  of the Soul}, he was already the author of an important book,
\emph{Novissima philosophia}, in which he attacked the Jesuit
theologians Oviedo, Arriaga, and Hurtado; {\ldots}%
\footnote{Armogathe and Carraud \cite[p.\;73]{Ar04}}
\end{enumerate}
Among the items censured by Spinula is the opinion that ``there is no
prime matter.''%
\footnote{Armogathe and Carraud \cite[p.\;75]{Ar04}.}

In the second half of the 17th century, several authors attempted to
reconcile atomism and catholicism,%
\footnote{Festa \cite[item~22]{Fe99}.}
including Donato Rossetti (1633--1686) \cite{Go96}.

In 1671, the Holy See alerted the archbishop of Naples concerning the
existence of attempts to reconcile the Cartesian philosophical system
and the eucharist, in the following terms:
\begin{enumerate}\item[]
``[I]l y avait dans cette ville [i.e., Naples] des gens qui, voulant
  prouver leur sup\'eriorit\'e, se faisaient les promoteurs des
  opinions philosophiques dun certain Renato des Cartes [i.e., Ren\'e
    Descartes] qui, il y quelques ann\'ees, a fait imprimer un nouveau
  syst\`eme philosophique, r\'e\-veillant ainsi les anciennes opinions
  des Grecs concernant les \emph{atomes}%
\footnote{Such a claim concerning Descartes is debatable; see
  note~\ref{f19}.}
{\ldots}; or, certains th\'eologiens pr\'etendent prouver, \`a partir
de cette doctrine, la mani\`ere dont les accidents du pain et du vin
se conservent apr\`es que ce pain et ce vin sont chang\'es en corps et
sang {\ldots}''%
\footnote{Translation by Festa \cite[item~28]{Fe99}; emphasis on
  ``atomes'' added.}
\end{enumerate}
A trial of atheists commenced in Naples in 1688 \cite[item~29]{Fe99};
see further in Fiorelli \cite{Fio21}.  A number of atomists came under
scrutiny, including the jurist Francesco D'Andrea (1625--1698) who was
acquitted nine years later.%
\footnote{Festa \cite[item~29]{Fe99}.}
%

%
%

\subsection{Sourcing opposition to atomism in Scotus}
\label{s47}

In 1651, indivisibles were placed on the jesuit list of permanently
banned doctrines~\cite[p.\;329]{He96}.
%
%
As noted by Leijenhorst and L\"uthy,
\begin{enumerate}\item[]
The incompatibility of atomism with the accepted interpretation of the
Eucharist was invoked, for example, {\ldots}~in the Jesuit
\emph{Ordination} [i.e., \emph{Ordinatio pro studiis superioribus}] of
1651, {\ldots}%
\footnote{Leijenhorst and L\"uthy \cite[p.\;396]{Le02}; cf.\;Vanzo
  \cite[p.\;211]{Va19}.}
%
%
%
\end{enumerate}
Leijenhorst and L\"uthy trace catholicism's opposition to atomism back
to Duns Scotus (c.\;1265/66--1308):
\begin{enumerate}\item[]
Ever since Duns Scotus had interpreted the decisions of the Fourth
Lateran Council of 1215 as entailing that the Church had formally
established transubstantiation as the canonical interpretation of the
Eucharist, physical atomism was generally viewed as a heresy.%
\footnote{Leijenhorst and L\"uthy \cite[p.\;396]{Le02}.}
\end{enumerate}
Cross notes that ``Scotus is fiercely opposed to any sort of atomism.''%
\footnote{Cross \cite[p.\,118]{Cr98}.  See Ariew \cite{Ar12} on
  disagreement between Aquinas and Scotus concerning the nature of
  prime matter (p.\,187) and on Descartes's anti-atomism (p.\,193).}
The final years of the pontificate of Innocent 3 \mbox{(1160/1--1216)}
witnessed a tightening of doctrinal control.  

The turn of the 12th century was a troubled time for catholicism,
which had to fend off a Manichaean challenge (the so-called
Albigensian heresy), against which Innocent 3 declared a crusade.
Besides the endorsement of transubstantiation as the canonical account
of the eucharist at the 4th Lateran Council in 1215 (see above), one
could mention a prohibition against the teaching of Aristotle's
\emph{Metaphysics} and commentaries by Avicenna and others,
promulgated in the same year by the papal legate Robert of Cour\c{c}on
in Paris.%
\footnote{De Vaux \cite[p.\;45]{De34}, Bertolacci
  \cite[p.\;214]{Be12}.}
Lahey notes that in the 14th century,
\begin{enumerate}\item[]
Wyclif seems to have employed his conception of spatio-temporal
indivisibles in arguing that the traditional doctrine of
transubstantiation was, as normally held, impossible, {\ldots}%
\footnote{Lahey \cite[p.\,102]{La09}.}
\end{enumerate}
Thus the tensions between indivisibles and transubstantiation were by
no means a new problem in the 17th century.

In sum, the objections to indivisibles formulated by jesuit
mathematicians (see Section~\ref{s9}) were only part of a greater
picture that involved acute doctrinal problems that had vexing
consequences for leading scholars such as Galileo and Descartes (and,
as we argue, the jesuate mathematicians) as well as a slew of
less prominent scientists.

\section{Clemency, Pasha, and cannons}
\label{s10}

In this section, we provide a non-mathematical explanation of the
suppression of the jesuates that obviates the need for Alexander's
unconvincing explanation analyzed in Section~\ref{s46}.

The established pattern is that, while the jesuits published a series
of decrees and books against indivisibles (see Sections~\ref{s9} and
\ref{s4}), the jesuates published a series of books on the geometry of
indivisibles: Cavalieri published two books (see Section~\ref{s2}) and
his student Degli Angeli published nine (see Section~\ref{s6}).  Both
of these scholars were highly placed in the jesuate order,%
\footnote{Thus, ``about 1652 [degli Angeli] was appointed prior of the
  monastery of the Gesuati in Venice, and shortly afterward he was
  given the post of provincial definer, a position he held until Pope
  Clement IX suppressed the order in 1668'' \cite[p.\,164]{Ca81}.}
a fact that surely colored their adversaries' attitude toward the
order itself.  The parties interested in making sure that the jesuate
problem ``penetrated the supreme ears'' (see Section~\ref{s51}) were
likely to include jesuits.

Cavalieri died in 1647, but Cavalieri's student Degli Angeli saw their
order of the jesuates suppressed with little clemency by pope
Clement\;9, on 6\;decem\-ber\;1668, by dint of a papal bull.%
\footnote{H\'elyot \cite[p.\;623]{He63}.  Alexander erroneously refers
  to this as a \emph{papal brief} in \cite[pp.\,171, 323]{Al15}.}
Alexander's hypothesis that the suppression of the jesuate order may
have been in reprisal to their members' work on indivisibles, is less
well established.%
\footnote{``[The jesuits] turned {\ldots}\;to the papal Curia in Rome,
  where their influence was decisive.  They could not punish [Degli]
  Angeli directly, so they let their fury rain on the order that
  sheltered him and his late teacher'' \cite[p.\,173]{Al15}.}

\subsection{Looming suppression}
\label{s51}

In 1644 pope Urban 8 died.  As reported by historian Uccelli, rumors
of a possible suppression started circulating no later than 1646, when
the General of the jesuates pleaded with the current pope in the
following terms:
\begin{enumerate}\item[]
From the beginning until these days there has not been a Pontiff who,
as a loving father, has not loved and helped it [i.e., the jesuate
  order]; for never have wrongdoings, that would discredit it,
penetrated those supreme ears.  Only today there are rumors of
removing our habit.%
\footnote{Uccelli \cite[p.\;90]{Uc65}.  Gagliardi attributes a similar
  passage to Cavalieri in \cite[p.\;483]{Ga04}.}
\end{enumerate}
Apparently the General was alarmed that damaging rumors concerning the
jesuates were beginning to penetrate the supreme ears.  Uccelli also
noted that a similar letter had been sent by Cavalieri from Bologna on
2 may 1646 to the Grand Duke Ferdinando II de' Medici (1621--1670) of
Tuscany, in recommendation of the same congregation.

The official historian of catholicism, Ludwig Pastor, glossed over the
suppression%
\footnote{``On August 18th, 1668, Venice was authorized to alienate
  Church property in aid of the Turkish war, and, subsequently, to
  raise a million ducats from suppressed monasteries'' (Pastor
  \cite[p.\;422]{Pa31}).  The ``monasteries'' in question remain
  unspecified, and the suppression of the jesuate order, unreported by
  Pastor.}
in his volume 31 where he had devoted over 100 pages to the papacy of
Clement\;9.%
\footnote{Pastor \cite[pp.\;314--430]{Pa31}.}
By contrast, his volume 1 does duly mention that in~1367,
\begin{enumerate}\item[]
Giovanni Colombini, the founder of the Gesuati, and his religious came
as far as Corneto to meet the Pope [Urban\;5], singing hymns of
praise.%
\footnote{Pastor \cite[p.\;95]{Pa1}.  What matters of course is not
that the jesuats were chanting hymns, but rather that the
\emph{founding} of the jesuates is covered in detail (11 lines) by
Pastor.  Significantly, the \emph{suppression} of the jesuates is not
even mentioned by Pastor, suggesting that it was perhaps a nettlesome
issue he preferred to avoid.}
\end{enumerate}
Pastor mentioned the founding of the jesuate order, but failed to
mention its suppression.  Indeed, the reader of the official
\emph{History of the popes} will find nothing to disabuse him of the
idea that Colombini's order still existed in 1940, publication date of
volume 31 of the English translation of Pastor's work.

\subsection{The \emph{aquavitae} issue}

There is some debate concerning the appellation \emph{Aquavitae
  Fathers} (or \emph{Brothers}).  The 1911 Encyclopedia Britannica
(EB) claims the following:
\begin{enumerate}\item[]
[In] the 17th century the Jesuati became so secularized that the
members were known as the Aquavitae Fathers, and the order was
dissolved by Clement IX in 1668.  \cite{En11}
\end{enumerate}
Thus the EB posits a causal link between alleged `secularisation' of
members and the appellation \emph{Aquavitae Fathers}, apparently
insinuating a charge of uncatholic dealings in alcohol.

%
%
Another account has it that in the hospitals, the jesuates distilled
liquors and gave them to the sick in order to help them bear the
pain~\cite{Ea11}.
%
%
Accordingly, the term \emph{Aquavitae Fathers} as applied to the
jesuates is linked to their deeds of charity.  As reported in the
Catholic Encyclopedia,
\begin{enumerate}\item[]
[I]n 1606 the Holy See allowed the reception of priests into the
congregation.  \cite{Ca13}
\end{enumerate}
Such a boost in their status is unlikely to have occurred had the EB's
linkage been correct.  Furthermore, in 1640, Urban\;8 approved new
constitutions for the jesuate order.%
\footnote{H\'elyot \cite[p.\;623]{He63}.}

The obscurity of the reasons sometimes given for the suppression of
the order (``some abuses crept in'')%
\footnote{\label{f11}Thus, Eves in his entry~$220^\circ$ starting with
  the words ``Written histories often contain hidden perpetuated
  errors'' \cite[p.\;35]{Ev69}, claimed that ``Certain abuses later
  made their way into the Jesuat order'' and proceeds to speculate
  darkly concerning ``the manufacture and sale of distilled liquors,
  apparently in a manner not acceptable to Canon Law,
  etc.''\;\cite[p.\;36]{Ev69}, without however providing any evidence
  as to his Canon Law claim.  Simmons lodged a similar allegation
  concerning liquors contrary to Canon Law \cite[p.\;106]{Si92}, still
  without evidence.  The canon involved, namely Council of Trent,
  Session 13, Canon 2 (see Section \ref{f7}) may have been at odds
  with the manufacture -- not of liquors but of theorems -- by leading
  jesuate scholars.}
suggests that they may have been introduced after 1668 so as
retroactively to justify the suppression -- coming as it did from an
infallible source.%
\footnote{Gagliardi claims that ``the Order was subsequently
suppressed by Clement IX in 1668 due to an insufficiency of vocations.
The suppressive records indeed reveal a substantial number of friars,
yet a meager count of novices and young aspirants''
\cite[p.\;36]{Ga20}.  It is open to question what the cause was and
what the effect: the jesuates having been the target of a negative
campaign for at least two decades prior to the ultimate suppression
(see Section~\ref{s51}), one can wonder whether such a campaign may
have had an effect on potential recruits.}

\subsection{Crusading for Candia}

Hippolyte (Pierre) H\'elyot (1660--1716) was eight years old when the
jesuates were suppressed.  Some decades later, H\'elyot detailed the
circumstances of the suppression as follows:
\begin{enumerate}\item[]
[D]ans l'Etat de Venise ils \'etaient assez riches, ce qui fit que la
r\'epublique demanda leur suppression \`a Cl\'ement IX, afin de
profiter de leurs biens, qui furent employ\'es \`a soutenir la guerre
que cette r\'epublique avait contre les Turcs qui assi\'egeaient pour
lors Candie; ce que le pape accorda l'an 1668, {\ldots}%
\footnote{H\'elyot \cite[p.\;623]{He63}.}
\end{enumerate}
The reference is to the Fifth Ottoman--Venetian War, or the Cretan War
(1645--1669); and to Candia, the capital of Crete at the time.
Uccelli provides a similar account.%
\footnote{Uccelli \cite[p.\;91]{Uc65}.}
Thus, Clement~9 appears to have determined that ``some abuses crept
in,'' suppressed the order, and transferred its assets to the Venetian
Republic,%
\footnote{Uccelli also notes the following: ``The Grand Duke
  [Ferdinando II de' Medici of Tuscany] wrote on December 18 of that
  year to his ambassador in Rome, so that he could interpose himself
  with his Holiness: so that `in the distribution of the effects and
  goods of similar monasteries existing in this dominion, S. B. may be
  to regard with the eyes of his supreme piety the needs of Tuscany,
  as he has regarded those, even if more urgent, of the Venetian
  dominion'{}'' \cite[p.\;91]{Uc65}.  The tenor of Ferdinando II's
  plea confirms that the papal bull consigned the jesuate assets to
  the Venetian republic even outside its realm.}
at the precise moment the latter was the most in need of funds to
fight K\"opr\"ul\"u Faz{\i}l Ahmed Pasha, the Ottoman Grand Vizier.

It is instructive to compare, following H\'elyot, two of the orders
suppressed by the 1668 papal bull: the \emph{Georges in Algha} and the
jesuates.  With regard to the former, H\'elyot writes:
\begin{enumerate}\item[]
Ils se sont bien \'eloign\'es dans la suite de la pauvret\'e et de
l'humilit\'e dont leurs fondateurs avaient fait profession et dont ils
leur avaient laiss\'e l'exemple.%
\footnote{H\'elyot \cite[p.\;400]{He63}.}
\end{enumerate}
H\'elyot goes on to detail their arrogant ways, and concludes:
\begin{enumerate}\item[]
[C]'est avec raison que Cl\'ement IX les supprima en 1668 et donna
tous leurs biens \`a la r\'epublique de Venise pour s'en servir dans
la guerre qu'elle avait contre les Turcs.%
\footnote{Ibid.}
\end{enumerate}
By contrast, H\'elyot details the rigors of the daily routine of the
jesuates, including prayers, dietary restrictions, and fasts,%
\footnote{H\'elyot \cite[p.\;623]{He63}.}
and, unlike the case of the \emph{Georges in Algha}, reports --
\emph{without endorsing~it} -- their 1668 suppression.  Possibly the
remarks in the EB are due to a confusion of the two orders.

\subsection{Danger of Venice making peace}

In a crusading effort, pope Clement~9 was able to mobilize contingents
from several European nations in 1668--69 to fight the infidels (the
sentiment was surely mutual).  Thus,
\begin{enumerate}\item[]
[Clement 9] met the wishes of [France and Spain with regard to the
  creation of Cardinals] on August 5th, 1669, in order to secure their
help for the Turkish war.%
\footnote{Pastor \cite[p.\;345]{Pa31}.}
\end{enumerate}
As a further illustration of papal savvy in contemporary
\emph{Realpolitik}, Pastor notes the following:
\begin{enumerate}\item[]
It was an act of sheer favour on the part of the Pope if he had taken
into consideration the wishes of France and Spain; {\ldots} it was
hoped that this act of grace would spur on Louis XIV to more intense
military action on behalf of Crete for \emph{there was danger of
  Venice making peace with the Turks}, thereby exposing both Italy and
the Emperor's own territories to Turkish aggression.%
\footnote{Pastor \cite[p.\;347]{Pa31}; emphasis added.}
\end{enumerate}
Pastor relates further that
\begin{enumerate}\item[]
Towards the end of May 1668, the papal galleys, under Vincenzo
Rospigliosi, put to sea in order to effect their junction with the
Venetian fleet commanded by Francesco Morosini.  Soon after Abbot
Airoldi was sent as special agent to obtain help for Crete from the
Italian and German princes.%
\footnote{Pastor \cite[p.\;419]{Pa31}.}
\end{enumerate}
Pastor furnishes a detailed account of all the Princes, Knights, and
Dukes participating in the war effort,%
\footnote{Pastor \cite[pp.\;420--421]{Pa31}.}
as well as a papal quote of \emph{Numbers}~10:35 in support of the
effort.%
\footnote{Pastor \cite[p.\;424]{Pa31}.}

Concerning the Ottoman side of the conflict in the 1660s, Baer writes:
\begin{enumerate}\item[]
Ottoman military policy, especially regarding the siege of Crete,
showed little sign of success, and the Ottomans, who considered
themselves to be the protectors of Islam, were demoralized.%
\footnote{Baer \cite[p.\,160]{Ba04}.}
\end{enumerate}
Such difficulties may have contributed to the Grand Vizier's eagerness
to conclude a peace treaty with Captain General Francesco Morosini.
The showdown between K\"opr\"ul\"u Faz{\i}l Ahmed Pasha and pope
Clement\;9 ended in an Ottoman victory by sword and cannon at Candia
on 6~september~1669,%
\footnote{Pastor \cite[p.\;427]{Pa31}.}
less than a year after the transfer of the jesuate assets to the
Venetians.

\section{Gregory's books held in great esteeme, but}
\label{s6}

Whatever the reasons for the summary suppression of the 300-year-old
jesuate order founded by the Blessed Giovanni Colombini,%
\footnote{Pope Gregory 13 included Colombini's name in the 1578
edition of the \emph{Martyrologium Romanum} (the official list of the
Saints, Blessed and martyrs).  The 1578 inclusion (over two centuries
after the recognition of the jesuates by Urban 5 and less than a
century before the suppression by Clement 9) endorsed a long-standing
practice in Siena (Colombini's birthtown), even though there was no
official act of beatification.}
collateral damage of the suppression included the books of James
Gregory (who visited Degli Angeli in the 1660s),%
\footnote{For details on Gregory's Italian trip see Crippa
  \cite{Cr20}.  Gregory visited Degli Angeli during the period
  immediately preceding the suppression of the jesuates.  The fact
  that, following the suppression of the jesuate order, Gregory's
  mathematical books were similarly suppressed in Italy, supports the
  hypothesis that the suppression of the jesuates may not have been
  unrelated to the mathematical investigations by their leading
  scholars.}
which were held ``in great esteeme, but not to be procured in Italy''
(Jean Bertet%
\footnote{Jean Bertet (1622--1692), jesuit, studied under Honor\'e
  Fabri, professor of mathematics at Aix \cite[p.\;506]{Wa14}.}
as reported by John Collins, quoted in Turnbull \cite[p.\,107]{Tu39};
see \cite{18f} for details).

Degli Angeli himself published nine books on the geometry of
indivisibles, including his \emph{De infinitis spiralibus inversis}
\cite{De67} which appeared merely a year before the suppression.%
\footnote{Giusti \cite[p.\;50 note\;39]{Gi80}.}
He did not write a word on them following the suppression, though he
remained active in mechanics, publishing a 79-page treatise in 1671
\cite{De71}, as well as in university teaching.

Decades later in the 1690s, we find Jacopo Riccati attending Degli
Angeli's course in astronomy at Padua.  In a touching episode of
historical continuity, the 72-year-old Degli Angeli encourages the
19-year-old Riccati to study the \emph{Principia Mathematica}.%
\footnote{Bertoloni Meli \cite[p.\;288]{Be06}.}

\section{Aftermath: \emph{le vide int\'egral} and a new dawn}
\label{s7}

Andersen notes that Degli Angeli, who was a jesuate like Cavalieri,
remarked in his \emph{De infinitis parabolis} \cite{De59} that the
circles opposed to the method of indivisibles mainly contained jesuit
mathematicians.%
\footnote{Andersen \cite[p.\;355]{An85}.}

In the aftermath of the hostility toward indivisibles, Italian
mathematics suffered a period of stunted growth in the emerging
infinitesimal analysis, and the centers of activity in what was
becoming mathematical analysis shifted north of the Alps.
Among the jesuits,
\begin{enumerate}\item[]
\ldots le grand nombre des math\'ematiciens {\ldots}~resta jusqu'\`a
la fin du XVIII$^e$ si\`ecle profond\'ement attach\'e aux m\'ethodes
euclidiennes.%
\footnote{Bosmans \cite[p.\;77]{Bo27}.}
\end{enumerate}
By 1700, Italy was a mathematical desert with regard to the new
mathematics:
\begin{enumerate}\item[]
en 1700, c'est le vide int\'egral en ce qui concerne la pratique des
math\'ematiques nouvelles en Italie\ldots{}%
\footnote{Robinet \cite[p.\,183]{Ro91}.}
\end{enumerate}

Less than four years after the suppression of the jesuates, Friedrich
Leibn\"utz's 26-year-old son arrived in Paris at the end of march
1672%
\footnote{See \cite[p.\,139]{An09}.}
on a diplomatic mission aimed at diverting Louis 14's appetite for
conquest away from the Dutch Republic and toward Egypt.  The mission
was aborted before it could start, the preparations for the French
military campaign being well underway%
\footnote{See \cite[pp.\,112--113]{Ly99}.}
by the time Gottfried Wilhelm reached Paris; but the arrival had
momentous significance in an entirely different realm.  Perhaps
sensing the doctrinal burden of the (overly) evocative term
\emph{indivisible}, Leibniz coined the term \emph{infinitesimal}%
\footnote{\label{f26}The term \emph{indivisibles} can be used in a
  narrower sense and in a broader sense.  Cavalieri's indivisibles
  possessed no thickness.  Torricelli started using indivisibles of
  variable thickness, i.e., already infinitesimals; see e.g.,
  \cite[pp.\;20--21]{An86}.  Many 17th century pioneers learned about
  Cavalieri's method through Torricelli; see further in \cite{Ba15}.
  Leibniz's sources were Blaise Pascal, John Wallis, and Honor\'e
  Fabri's \emph{Synopsis Geometrica} (1669).  It is not surprising
  that Leibniz sometimes referred to his own method as
  \emph{indivisibles} long after he introduced the new term
  \emph{infinitesimal}.  Adhering strictly to the narrower sense of
  \emph{indivisibles} (as codimension~$1$ entities) is an innovation
  of 20th century scholars starting with Koyr\'e.  It makes for a neat
  classification but it does not faithfully reflect the variety of
  17th century practice.}
in 1673, taking up a proposal by Nicolaus Mercator.%
\footnote{See Probst \cite[p.\;200]{Pr18}, who also notes that Wallis
  used the term \emph{pars infinitesima} already in 1670.}
Within a short few years, infinitesimal calculus was born.

Leibnizian infinitesimals were controversial in their own right, both
in Leibniz's time and in current Leibniz scholarship.%
\footnote{The related bibliography is vast.  One could mention Bos's
  article \cite{Bo74} and the articles cited in note~\ref{f143}.}
Did the jesuits' sacramentally inspired opposition to atomism and
indivisibles also apply to the new technique of infinitesimals?
Sherry argues that it did not necessarily:
\begin{enumerate}\item[]
indivisibles but not infinitesimals conflict with the doctrine of the
Eucharist, the central dogma of the Church.%
\footnote{Sherry \cite[p.\;367]{18k}.}
\end{enumerate}
A vast correspondence between jesuit des Bosses and Leibniz
\cite{Le07} attests to the jesuit's keen interest in Leibniz's work.%
\footnote{\label{f143}Bair et al.~\cite[Section 4.10]{18a}.
Interpreting Leibnizian infinitesimals is an area of lively debate.
In 2021, Bair et al.~published a comparative study of three
interpretations \cite{21a}; Katz et al.~\cite{21g} presented three
case studies in Leibniz scholarship.  In 2022, Katz et al.~presented
and analyzed a pair of rival approaches \cite{22b}; Archibald et
al.~formulated some criticisms \cite{Ar22b}.  Bair et al.~provided a
brief response \cite{23a} and a detailed response \cite{22a}.  Recent
work includes \cite{23g}, \cite{24b} and \cite{24c}.  Arthur and
Rabouin formulated further objections in \cite{Ar24}.  In their latest
book \cite{Ar25}, they double down on their interpretation following
Ishiguro \cite[ch.\;5]{Is90}, and claim that Leibnizian infinitesimals
were not \emph{quantities}.  A rebuttal appears in Katz and Kuhlemann
\cite{23f}.}
Indivisibles were also taught by Cavalieri's successor Pietro Mengoli
\cite{Me59}, as well as by the jesuits Fabri \cite{Fa69} and Dechales
\cite{De74}.

\section{Final remarks}
\label{s8}

The theological underpinnings of the 17th century debate over
indivisibles require further study.  It may be of interest to analyze
the philosophical/theological pedigree of the jesuates in relation to
their being more accepting of the emerging mathematics of
indivisibles.

The indivisibles were controversial within the jesuit order.  While
some maverick jesuits, including Gregoire de Saint Vincent, attempted
to exploit them (Section~\ref{s9}), mathematicians like Paul Guldin,
Andr\'e Tacquet, and others opposed them (Section~\ref{s9}), and their
doctrinal superiors in the order issued a series of bans aimed against
teaching indivisibles in jesuit colleges (Section~\ref{s41}).
Meanwhile, scholars in the jesuate order were among the leading
practitioners of indivisibles (Sections~\ref{s21} and \ref{s6}):
Cavalieri and his student Degli Angeli published a series of books on
indivisibles.  Cavalieri attempted to avert the looming suppression of
their order (Section~\ref{s51}).

Such a difference in attitude towards indivisibles between the two
orders may be related to differences in the theological pedigree of
the rival orders of the jesuates and the jesuits, respectively via
Augustine of Hippo (354--430) and Thomas Aquinas (1225--1274), back to
respectively Neoplatonists and Aristotle.  The relevance of such
pedigree is due to the connection of transubstantiation to what many
catholic theologians conceived of as Aristotelian doctrine of matter
and form, the latter being thought of, in Thomist theology, as
clashing with atoms and indivisibles.  The counter-reformation context
of the tensions over transubstantiation was analyzed by Festa
\cite{Fe90}, Hellyer \cite {He96}, McCue \cite{Mc68}, and others.  The
jesuit Pallavicino is on record concerning the tension between atomism
and the catholic interpretation of the eucharist (Section~\ref{f7}),
while jesuits Grassi and Inchofer attacked Galileo with the weapon of
eucharistic theology (Section~\ref{s43}).

H\'elyot traces the theological pedigree of the jesuates back to ``la
r\`egle de saint Augustin'' \cite[p.\;620]{He63}.
%
%
%
Some of the Italian atomist catholic scholars were influenced by
Descartes, whose doctines were deemed to be at odds with Thomist
theology by the catholic hierarchy (Section~\ref{s45}).  A comparison
of the intellectual roots of the jesuates and the jesuits in the
thought of respectively Augustine and Thomas merits further analysis.

\section*{Acknowledgments} We are grateful to Roger Ariew, Ugo
Baldini, Davide Crippa, Vladimir Kanovei, Paolo Palmieri, Siegmund
Probst, and Katerina Trlifajova for subtle suggestions.

\end{document}